\documentclass{article}
\usepackage{amsmath, epsfig,amssymb, multicol}
\usepackage[active]{srcltx}
\usepackage{color}
\newtheorem{lemma}{Lemma}
\newtheorem{theorem}{Theorem}

\renewcommand{\Pr}{{\rm Pr}}
\newcommand{\E}{{\rm E}}
\newcommand{\C}{{\cal C}}

\newcommand{\F}{{\cal F}}
\newcommand{\h}{{\cal H}}
\renewcommand{\P}{{\cal P}}
\newcommand{\Q}{{\cal Q}}

\newcommand{\W}{{\cal W}}
\newcommand{\X}{{\cal X}}
\renewcommand{\exp}{{\rm exp}}
\newcommand{\bp}{{\rm bp}}
\newcommand{\llll}{\log\log \log \log}
\renewcommand{\lll}{\log\log\log}

\title{On the decomposition of random  hypergraphs}
\author{Xing Peng
\thanks{Center for Applied Mathematics, Tianjin University, Tianjin, 300072, P.~R.~China,
({\tt x2peng@tju.edu.cn}).}}
\date{}
\begin{document}
\maketitle

\begin{abstract}
For an $r$-uniform hypergraph $H$, let $f(H)$ be the minimum number of  complete $r$-partite $r$-uniform subhypergraphs of $H$ whose edge sets partition the edge set of $H$. For a graph $G$, $f(G)$ is the bipartition number of $G$ which was introduced by Graham and Pollak in 1971. In 1988, Erd\H{o}s conjectured that  if $G \in G(n,1/2)$, then with high probability $f(G)=n-\alpha(G)$, where $\alpha(G)$ is the independence number of $G$. This conjecture and related problems have received a lot of attention recently.
In this paper, we study the  value of $f(H)$ for a typical $r$-uniform hypergraph $H$. More precisely, we prove that if $(\log n)^{2.001}/n \leq p \leq 1/2$ and $H \in H^{(r)}(n,p)$, then with high probability $f(H)=(1-\pi(K^{(r-1)}_r)+o(1))\binom{n}{r-1}$, where $\pi(K^{(r-1)}_r)$ is the Tur\'an density of $K^{(r-1)}_r$.
\end{abstract}

\section{Introduction}
For a graph $G$, the {\it bipartition number} $\tau(G)$ is the minimum number of complete bipartite subgraphs of $G$ so that each edge of $G$ belongs to exactly one of them. This parameter of a graph was introduced by Graham and Pollak \cite{gp1} in 1971. The famous  Graham--Pollak \cite{gp1} Theorem asserts $\tau(K_n)=n-1$. Since it's original proof using Sylvester's Law of Intertia, many other proofs have been discovered, see \cite{lo}, \cite{peck}, \cite{tv}, \cite{v1}, \cite{v2}, \cite{yan}.

Let $\alpha(G)$ be the independence number of $G$. It is easy to observe $\tau(G) \leq |V(G)|-\alpha(G)$.
Erd\H{o}s (see \cite{krw}) conjectured that the equality holds for almost all graphs. Namely,  if $G \in  G(n,1/2)$, then  $\tau(G)=n-\alpha(G)$  with high probability.  Alon \cite{alon2} disproved this conjecture by showing $\tau(G) \leq n-\alpha(G)-1$  with high probability   for most values of $n$. Alon's upper bound on the bipartition number of random graphs $G \in G(n,1/2)$ was improved by
Alon, Bohman, and Huang \cite{alon3} recently.  Chung and the author proved that if $G \in G(n,p)$,  $p$ is a constant, and $p \leq 1/2$, then with high probability we have $\tau(G) \geq n- \delta(\log_{1/p} n)^{3+\epsilon}$ for any constants $\delta$ and $\epsilon$. When $p$ satisfies $\tfrac{2}{n} \leq p \leq c$ for some absolute (small) constant $c$, Alon \cite{alon1} showed that if $G \in G(n,p)$, then $\tau(G)=n-\Theta\left( \tfrac{\log(np)}{p}\right)$ with high probability.

The hypergraph analogue of the bipartition number is well-defined.
For $r\geq 3$ and an  $r$-uniform hypergraph $H$, let $f(H)$ be the minimum number of   complete $r$-partite $r$-uniform subhypergraphs of $H$ whose edge sets partition the edge set of $H$.
Aharoni and Linial (see \cite{alon1}) first asked to determine the value of $f(K_n^{(r)})$ for $r \geq 3$, where $K_n^{(r)}$ is the complete $r$-uniform hypergraph with $n$ vertices.  The value of $f(K_n^{(r)})$  is related to a perfect hashing problem from computer science.
 Alon \cite{alon1} proved $f(K^{(3)}_n)=n-2$ and  $c_1(r) n^{\lfloor \tfrac{r}{2} \rfloor} \leq f(K^{(r)}_n) \leq c_2(r) n^{\lfloor \tfrac{r}{2} \rfloor}$ for $r \geq 4$.
For  improvements and variations, readers are referred to \cite{ck}, \cite{cktv}, \cite{ckv}, and \cite{ct}.
In this paper, we examine the  value of $f(H)$ for the random hypergraph $H \in H^{(r)}(n,p)$. To state our main theorem, we need few more definitions.

For an $r$-uniform hypergraph $H$, the {\it Tur\'an number} $\textrm{ex}(n,H)$ is the maximum number of edges in an $n$-vertex $r$-uniform hypergraph which does not contain $H$ as a subhypergraph.  We define the {\it Tur\'an density} of $H$ as
\[
\pi(H)=\lim_{n \to \infty} \frac{\textrm{ex}(n,H)}{\binom{n}{r}}.
\]
For each $r \geq 3$,  we use  $K^{(r-1)}_{r}$ to denote the compete $(r-1)$-uniform hypergraph with $r$ vertices.

By extending techniques from \cite{alon2} and \cite{cp}, we are able to prove the following  theorem.

\begin{theorem} \label{main}
For  $r\geq 3$, if $(\log n)^{2.001}/n \leq p \leq 1/2$  and $H \in H^{(r)}(n,p)$,   then  with high probability we have
\[
f(H)=(1-\pi(K^{(r-1)}_r)+o(1)) \binom{n}{r-1}.
\]
\end{theorem}
From this theorem, we can see the typical value of $f(H)$  has the  order of magnitude $n^{r-1}$  while $f(K^{(r)}_n)$ has the order of magnitude $n^{{\lfloor \tfrac{r}{2} \rfloor}}$.  We note $\pi(K_3^{(2)})=\tfrac{1}{2}$ while  the value of $\pi(K^{(r-1)}_r)$ is not known for $r \geq 4$.
We remark here that our techniques also work for  $p \leq 1-c$ for any small positive constant $c$. However, we restrict out attention to the case where $p \leq 1/2$ in this paper.

We will use the following notation throughout this paper.
  For each $r \geq 3$, we will use $[n]$ to denote the set $\{1,2,\ldots,n\}$ and $\binom{[n]}{r}$ to denote the collection of all $r$-subsets of $[n]$. If $A_1,A_2,\ldots,A_r$ are pairwise disjoint subsets of $[n]$, then we use $\prod_{i=1}^r A_i$ to denote those $r$-subsets $F$ of $[n]$ such that $|F \cap A_i|=1$ for each $1 \leq i \leq r$. We may also write $A_1 \times A_2 \times \cdots \times A_r$ for $\prod_{i=1}^r A_i$ on some occasions. The complete $r$-partite $r$-uniform hypergraph whose  vertex parts are  $A_1,A_2,\ldots,A_r$ is the $r$-uniform hypergraph with the edge set $\prod_{i=1}^r A_i$.

  Let $H$ be an $r$-uniform hypergraph with vertex set $[n]$ and edge set $E$.
For pairwise disjoint subsets $A_1,A_2,\ldots,A_r \subset [n]$, we say $A_1,A_2,\ldots,A_r$ form a complete $r$-partite $r$-uniform hypergraph if $\prod_{i=1}^r A_i \subseteq E(H)$.

For an $r$-uniform hypergraph $H$, suppose $E(H)=\sqcup_{i=1}^{q} \prod_{j=1}^r A_j^i$ is a partition of the edge set of $H$. For each $1 \leq i \leq q$, the $i$-th  complete  $r$-partite $r$-uniform hypergraph $H_i$  has vertex parts $A_1^i,\ldots,A_r^i$. We always assume $|A_1^i|\leq \cdots \leq |A_r^i|$. We say $H_i$ is a trivial complete $r$-partite $r$-uniform hypergraph if $|A_1^i|=\cdots=|A_{r-1}^i|=1$. Otherwise, we say $H_i$ is a nontrivial one.
  The {\it prefix} $P_i$ of $H_i$ is the set  $\{A_1^i,\ldots,A_{r-1}^i\}$ and the prefix set $\P$ of the partition is $\{P_1,\ldots,P_q\}$.

  We will use $H^{(r)}(n,p)$  to denote the random $r$-uniform hypergraph in which each $r$-set $F \in \binom{[n]}{r}$ is selected as an edge with probability $p$ independently.
 We say an event $\X_n$  occurs with high probability  if  the probability that $\X_n$ holds goes to one as $n$ approaches infinity. All logarithms are in base 2, unless otherwise specified.

 The outline of the proof is the following. We
 will prove  $f(H) \leq (1-\pi(K^{r-1}_r)-\epsilon)\binom{n}{r-1}$  with small probability  for any  positive  constant $\epsilon$. To do so, for a given prefix set $\P=\{P_1,\ldots,P_q\}$ with $q \leq (1-\pi(K^{r-1}_r)-\epsilon)\binom{n}{r-1}$,  let $\P_1=\{P_i \in \P: |P_j^i|=1 \text{  for each } 1 \leq j \leq r-1\}$ and $\P_2=\P \setminus \P_1$.
    We will show that there are at least $c(\epsilon)n^r$ edges of  $H \in H^{(r)}(n,p)$ which must be covered by some nontrivial complete $r$-partite $r$-uniform hypergraph with the prefix from $\P_2$.  Theorem \ref{t:submain} will tell us this probability is sufficiently small.   We will prove an upper bound for the number of possible choices for $\P$ and apply the union bound to complete the proof.

 The rest of the paper is organized as follows. In Section 2, we will prove several necessary lemmas. In Section 3, we will present the proof of   an auxiliary theorem which is the key ingredient in the proof of the main result. Theorem \ref{main} will be proved in Section 4. Few concluding remarks will be mentioned in Section 5.

\section{Lemmas}
In this section, we will collect some necessary lemmas which are needed to prove the main theorem.
We will  use the following versions of Chernoff's inequality and   Azuma's  inequality.
\begin{theorem}{\cite{chernoff}}
\label{t:chernoff}
 Let $X_1,\ldots,X_n$ be independent random variables with
$$\Pr(X_i=1)=p_i, \qquad \Pr(X_i=0)=1-p_i.$$
We consider the sum $X=\sum_{i=1}^n X_i$
with expectation $\E(X)=\sum_{i=1}^n p_i$. Then we have
\begin{eqnarray*}
\mbox{(Lower tail)~~~~~~~~~~~~~~~~~}
\qquad \qquad  \Pr(X \leq \E(X)-\lambda)&\leq& e^{-\lambda^2/2\E(X)},\\
\mbox{(Upper tail)~~~~~~~~~~~~~~~~~}
\qquad \qquad
\Pr(X \geq \E(X)+\lambda)&\leq& e^{-\frac{\lambda^2}{2(\E(X) + \lambda/3)}}.
\end{eqnarray*}
\end{theorem}

\begin{theorem} \cite{azuma}
 \label{t:azuma}
Let $X$ be a random variable determined by $m$ trials $T_1,\ldots,T_m$, such that for each $i$, and any two possible sequences of outcomes $t_1,\ldots,t_{i-1},t_i$ and $t_1,\ldots,t_{i-1}, t_i'$:
\[
|\E\left(X|T_1=t_1,\ldots,T_i=t_i\right) -\E\left(X|T_1=t_1,\ldots, T_{i-1}=t_{i-1},T_i=t_i'\right)   | \leq c_i
\]
then
\[
\Pr\left(|X-\E(X)| \geq \lambda  \right) \leq 2 \exp\left(-{\lambda}^2/2\sum_{i=1}^m c_i^2\right).
\]
\end{theorem}


Recall that if $A_1,\ldots,A_r$ form a complete $r$-partite $r$-uniform hypergraph, then we  assume $|A_1| \leq |A_2| \leq \cdots \leq |A_r|$.  We have the following lemma.

\begin{lemma} \label{l:lm1}
For  $H \in H^{(r)}(n,p)$ with $p \leq 1/2$,
with high probability  the vertex parts $A_1,A_2, \cdots,   A_r$ of each complete $r$-partite $r$-uniform hypergraphs in $H$   satisfy  $\prod_{i=1}^{r-1} |A_i| < (r+1) \log n$.
\end{lemma}

\noindent
{\bf Proof:} We need only to prove the lemma for $p=1/2$. For  a collection of pairwise disjoint sets $A_1, A_2, \ldots, A_r \subset [n]$,   we assume $|A_i|=k_i$ for each $1 \leq i \leq r$ and $k_1 \leq k_2 \leq \cdots \leq k_r$.    Fix a selection of $A_1, \ldots, A_r$,
the probability that they form a complete $r$-partite  $r$-uniform hypergraph in $H^{(r)}(n,1/2)$  is
$2^{-\prod_{i=1}^r k_i}.$ For fixed $k_1,\ldots,k_r$,   there are at most $\prod_{i=1}^r \binom{n}{k_i}$
 choices for $A_1,A_2,\ldots,A_r$ such that $|A_i|=k_i$ for each $1 \leq i \leq r$.  Therefore, for fixed $k_1,\ldots,k_r$ satisfying $ \prod_{i=1}^{r-1}  k_i \geq (r+1) \log n$ and $k_1 \leq \cdots \leq k_r$, the probability that there are pairwise disjoint sets   $A_1,A_2,\ldots,A_r$  such that $|A_i|=k_i$ and they  form a complete $r$-partite  $r$-uniform hypergraph is at most
\begin{align*}
\prod_{i=1}^r \binom{n}{k_i} 2^{-\prod_{i=1}^r k_i} &<  2^{(\sum_{i=1}^r k_i) \log n -\prod_{i=1}^r k_i}\\
                                                           &=2^{k_r\left((\sum_{i=1}^{r-1}k_i/k_r+1)\log n-\prod_{i=1}^{r-1} k_i\right)}  \\
                                                           & \leq   2^{k_r(  r \log n-\prod_{i=1}^{r-1} k_i)}\\
                                                           &<  2^{-k_r \log n}
\end{align*}
Put   $s=\prod_{i=1}^{r-1} k_i$. We next estimate how many choices of $k_1,\ldots,k_r$ such that $\prod_{i=1}^{r-1} k_i=s$ and $k_1 \leq \cdots \leq k_r$.
Let $t=\sum_{i=1}^{r-1} k_i$.  If $s \geq \log n $,  then  $t  \leq s+r< 2s$ and $k_r \geq k_{r-1} \geq s^{1/{(r-1)}}$.  Thus the number of choices for $k_1,\ldots,k_{r-1}$ satisfying  $\prod_{i=1}^{r-1} k_i=s$ and $k_1 \leq \cdots \leq k_{r-1}$ is less than the number of positive solutions to the equation
 $\sum_{i=1}^{r-1} k_i=t$, which is   less than  $\binom{2s}{r-2}$    as  $t \leq 2s$. We have at most $n$ choices for $k_r$ regardless the choices of $k_1,\ldots,k_{r-1}$.
Therefore,  the probability that there are   $A_1,A_2,\ldots,A_r$  which satisfy $s=\prod_{i=1}^{r-1} |A_i| \geq (r+1) \log n$ and  form a complete $r$-partite  $r$-uniform hypergraph in $H^{(r)}(n,1/2)$ is at most
\[
\sum_{s=(r+1) \log n} ^n   n \binom{2s}{r-2}  2^{-k_r \log n} \leq  \sum_{s=(r+1) \log n} ^n  n  \binom{2s}{r-2} 2^{-s^{1/(r-1)} \log n} =o(1)
\]
  Then the lemma  follows from Markov's inequality.
 \hfill $\square$


For an $r$-uniform  hypergraph $H=(V,E)$  and a prefix $P=\{A_1, A_2, \ldots, A_{r-1}\}$, we define
 \[
 V(H,P)=\{v:  v \in V(H) \setminus (  \cup_{i=1}^{r-1} A_i)  \textrm{ and } F \in E(H) \textrm{ for each } F \in A_1 \times \cdots \times A_{r-1} \times \{v\} \}.
 \]
 \begin{figure}[ht]
 \centerline{\psfig{figure=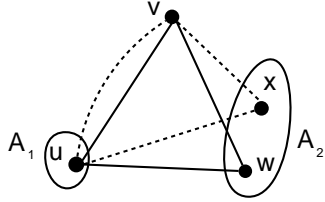, width=0.3\textwidth} }
\caption{An example with $r=3$, $P=\{A_1,A_2\}$, and $v \in V(H,P)$.}
\label{fg:fg1}
\end{figure}
Figure \ref{fg:fg1} is an illustrative example for $v \in V(H,P)$.  It  follows    that  $A_1, A_2, \ldots,A_r$ form a complete $r$-partite $r$-uniform hypergraph if $A_r$ is contained in $V(H,P)$, namely,  $A_r \subseteq V(H,P)$.
 We say an edge $F \in E(H)$ is covered by a complete $r$-partite $r$-uniform hypergraph with the prefix $P$ if $F \in A_1 \times \cdots \times A_{r-1} \times V(H,P)$.

Let $\P= \{P_1,\ldots,P_q\}$ be a prefix set, where $P_i=\{A_1^i,\ldots,A_{r-1}^i\}$.  We define $g(H,\P)$ as the set of edges of $H$ which are contained by exactly one  complete $r$-partite $r$-uniform hypergraphs whose  prefix is  from $\P$. It is easy to see
\[
g(H,\P) \leq \sum_{i=1}^q g(H,P_i)=\sum_{i=1}^q |V(H,P_i)| \prod_{j=1}^{r-1} |A_i^j|.
\]
We have the following lemma on  $g(H,\P)$ for $H \in H^{(r)}(n,p)$.
\begin{lemma} \label{l:lm2}
  Assume $p \leq 1/2$ and  $H \in H^{(r)}(n,p)$. Let $c(n)$ be a fixed function.  Given a prefix set  $\P=\{P_1,\ldots,P_q\}$, where $P_i=\{A_1^i,A_2^i,\ldots,A_{r-1}^i\}$  and $c(n)  \leq  \prod_{j=1}^{r-1} |A_j^i| < (r+1)\log n$ for each $1 \leq i \leq q$,
then we have
\[
\Pr\left( g(H,\P) \geq q c(n) p^{c(n)} n + 2n^{r-0.3} \right) \leq 2 \exp(-n^{r-0.8}).
\]

\end{lemma}

\noindent
{\bf Proof:} We shall use Theorem \ref{t:azuma} to prove this lemma.  Let $m=\binom{n}{r}$ and we list all $r$-sets of $[n]$ as $F_1,F_2,\ldots,F_m$. For each  $1 \leq i \leq m$, we consider $T_i \in \{\textrm{H}, \textrm{T}\}$, here $T_i=\textrm{H}$ means $F_i$ is an edge and  $T_i=\textrm{T}$ means $F_i$ is a  nonedge.
To simplify the notation, we use $X$ to denote the random variable $g(H,\P)$. We observe that $X$ is determined by $T_1,\ldots,T_m$. Fix the outcome $t_j$ of $T_j$ for each $1 \leq j \leq i-1$, we wish to show an upper bound for
\begin{equation} \label{eq:lm1}
\left|\E(X|T_1=t_1,\ldots,T_{i-1}=t_{i-1}, T_i=\textrm{H})-\E(X|T_1=t_1,\ldots,T_{i-1}=t_{i-1},T_i=\textrm{T})\right|
\end{equation}
If $T_i=\textrm{H}$, then we can assume $F_i$ is contained by some hypergraph whose prefix is $P_i$ for some $1 \leq i \leq q$. Otherwise the change of the outcome of $T_i$ will not effect the value of $X$. Suppose $F_i=\{v_1,\ldots,v_r\}$, where $A_j^i \cap F_i=\{v_j\}$ for each $1 \leq j \leq r-1$ and $v_r \not \in \cup_{j=1}^{r-1} A_j^i$. We next examine other edges which get covered because  we change $F_i$ as an edge. These edges are from the family  $A_1^i \times \cdots A_{r-1}^i \times \{v_r\}$. Therefore $\prod_{j=1}^{r-1} |A_j^i|$ is an upper bound for \eqref{eq:lm1}. Recalling the assumption $\prod_{j=1}^{r-1} |A_j^i| < (r+1)\log n$, then \eqref{eq:lm1} can be bounded from above by $(r+1)\log n$.

We note  $\E(g(H,P_i)) \leq c(n)p^{c(n)} n$ as we assume $\prod_{j=1}^{r-1} |A_j^i| \geq c(n)$. We get
\[
\E(X) \leq \sum_{i=1}^q \E(H,P_i) \leq c(n)qp^{c(n)}n.
\]
Applying Theorem \ref{t:azuma} with $\lambda=2n^{r-0.3}$ and $c_i=(r+1) \log n$, we obtain
\begin{align*}
\Pr\left(X \geq c(n)qp^{c(n)}n+2n^{r-0.3}\right) &\leq \Pr\left(X \geq \E(X) +2n^{r-0.3}\right)\\
                                    &\leq 2\exp\left(-4n^{2r-0.6}/(2m(r+1)\log n)\right) \\
                                     &\leq 2\exp(-n^{r-0.8})
\end{align*}
here we used the fact $m<n^r$. \hfill $\square$

We need a theorem which provides a lower bound for the number of uncovered edges.
Let $k(n)$ and $l(n)$ be  given functions. Suppose $\F  \subset \binom{[n]}{r}$ and $\Q$
be the power set of $\binom{[n]}{r} \setminus \F$. Consider a function $\C: \F \to \Q$ such that for each $F \in \F$ and each $R \in \C(F)$, we have $|R \cap F|=r-1$.   Let $h(H,\F,\C)$ be the number of $F \in \F$ such that $F$ is an edge in $H \in H^{(r)}(n,p)$ and $R$ is not an edge in $H \in H^{(r)}(n,p)$ for all $R \in   \C(F) $. We have the following lemma.

\begin{lemma} \label{l:lm3}
 Suppose $p \leq 1/2$ and $ \F \subset \binom{[n]}{r}$.  Assume $H \in H^{(r)}(n,p)$, $|\C(F)| \leq k(n)$ for each $F \in \F$,  and for each $R \in  \cup_{F \in \F} \C(F) $, the number of $F \in \F$ satisfying $R \in \C(F)$ is at most $l(n)$, here $l(n)$ and $k(n)$ are some given functions. Then we have
\[
\Pr\left(h(H,\F,\C) \leq  |\F|p(1-p)^{k(n)} -  2n^{r-0.01} \right) \leq 2\exp(-n^{r-0.02}/l(n)^2).
\]
\end{lemma}

\noindent
{\bf Proof:} To simplify notation, we use $X$ to denote the random variable $h(H,\F,\C)$ again. We list all $r$-sets from
 $\F \cup_{F \in \F} \C(F)$
 as $F_1,F_2,\ldots,F_m$, here $m \leq \binom{n}{r}$. For each $F_i$, we consider $T_i \in \{\textrm{H},\textrm{T}\}$, here $T_i= \textrm{T}$ means $F_i$ is an edge and $T_i=\textrm{T}$ means $F_i$ is not an edge. Given the outcome $t_j$ of $T_j$ for each $1 \leq j \leq i-1$ we wish to establish an upper bound for
\begin{equation} \label{eq:lb}
\left|\E(X|T_1=t_1,\ldots,T_{i-1}=t_{i-1}, T_i=\textrm{H})-\E(X|T_1=t_1,\ldots,T_{i-1}=t_{i-1},T_i=\textrm{T})\right|.
\end{equation}
If $F_i \in \F$, then changing the outcome of $T_i$ can only effect \eqref{eq:lb} by one. If $F_i \in \cup_{F \in \F} \C(F)$, then changing the outcome of $T_i$ can  effect \eqref{eq:lb} by at most $l(n)$ since $F_i \in \C(F)$ for at most $l(n)$  $r$-set $F$. Therefore, the expression \eqref{eq:lb} can be bounded from above by $l(n)$.  Applying Theorem \ref{t:azuma} with $\lambda=2n^{r-0.01}$ and $c_i=l(n)$, we  get
\[
\Pr\left( |X-\E(X)| \geq 2n^{r-0.01} \right) \leq 2\exp\left(-4n^{2r-0.02}/2\sum_{i=1}^m c_i^2\right)  \leq 2\exp(-n^{r-0.02}/  l(n)^2),
\]
here we used $m  \leq \binom{n}{r}$. We note $\E(X)= \sum_{F \in \F} p(1-p)^{|\C(F)|} \geq |\F| p (1-p)^{k(n)}$ as $|\C(F)| \leq k(n)$. Therefore
\begin{align*}
\Pr\left(h(H,\F,\C) \leq  |\F|p(1-p)^{k(n)} -  2n^{r-0.01}  \right) & \leq \Pr\left( |X-\E(X)| \geq 2n^{r-0.01} \right)\\
                                                                           &\leq 2\exp(-n^{r-0.02}/ l(n)^2)
\end{align*}
We proved the lemma. \hfill $\square$

When $p \leq 1/\llll n$,  we adapt the approach in  \cite{alon2}.  The following two lemmas are the hypergraph version of Lemma 3.1 and Lemma 3.2 in \cite{alon2}. Before we state them, we need one additional definition.

For  positive integers $m \geq \log n$ and $r \geq 3$, let ${\cal T}_m$ be the set of tuples $(a_1,a_2,\ldots,a_r)$ satisfying the following properties.
\begin{enumerate}
\item[1:] $a_i$ is a positive integer for each $1 \leq i \leq r$;
\item[2:] $1 \leq a_1 \leq a_2 \leq \cdots \leq a_r$;
\item[3:] $a_1 \cdots a_r=m$;
\item[4:]  $a_{r-1} \geq 2$.
\end{enumerate}

\begin{lemma} \label{l:lm8}
For any constant $c$,  if $p$ satisfies $(\log n)^{2.001}/n \leq p \leq 1/\lll n$, then  the following holds for  $n$ large enough. For every integer $m$ satisfying
\[
\frac{ pcn}{16} \leq m \leq \frac{pcn}{4},
\]
we have
\[
 \sum_{(a_1,\ldots,a_r) \in {\cal T}_m} \binom{n}{a_1} \binom{n-a_1}{a_2} \cdots \binom{n-\sum_{i=1}^{r-1}a_i}{a_r} p^m \leq 2^{-0.3 \log(1/p)m}.
\]
\end{lemma}
Recall that a complete $r$-partite $r$-uniform hypergraph whose vertex parts $A_1,\ldots,A_r$ satisfying $|A_1|\leq |A_2| \leq \cdots \leq |A_r|$  is nontrivial if $\prod_{i=1}^{r-1} |A_i| \geq 2$.

\begin{lemma}\label{l:lm9}
  For any constant $c$, if $p$ satisfies $(\log n)^{2.001}/n \leq p \leq 1/\lll n$, then the probability that $H \in H^{r}(n,p)$ contains a set of at most $2n^{r-1}$ nontrivial complete $r$-partite $r$-uniform hypergraphs which cover at least $pcn^r/4$ edges is at most $2^{-0.05pc \log(1/p)n^r}$.
\end{lemma}
As proofs of the two lemmas above go the same lines as those for proving Lemma 3.1 and Lemma 3.2 in \cite{alon2}, they are  omitted  here.

\section{An auxiliary theorem}
Let $\F \subset \binom{[n]}{r}$ with $|\F| \geq cn^r$ for some positive constant $c$.  Suppose the probability $p$ satisfies $1/\llll n \leq p \leq 1/2$.
We shall prove that if $H \in H^{(r)}(n,p)$, then with small probability that there are few nontrivial complete $r$-partite $r$-uniform hypergraphs such that each edge $F \in E(H) \cap \F$ is in exactly one of them.


\begin{theorem} \label{t:submain}
Assume $\F \subset \binom{[n]}{r}$ with $|\F| \geq c n^r$ for some positive constant $c$. Let $\P=\{P_1,\ldots,P_t\}$ be a given prefix set, where
 $t=|\P| \leq n^{r-1}$ and $P_i=\{A_1^i,\ldots,A_{r-1}^i\}$ satisfying $2 \leq \prod_{j=1}^{r-1} |A_j^{i}| <  (r+1) \log n$ for each $1 \leq i \leq t$.
If $1/\llll n \leq p \leq 1/2$ and  $H \in H^{(r)}(n,p)$, then with probability at most  $3\exp(-n^{r-0.92})$  there are $t$ nontrivial complete $r$-partite $r$-uniform hypergraphs  such that its  prefix set is $\P$ and each edge $F \in E(H) \cap \F$ is in exactly one of these hypergraphs.
 \end{theorem}


Suppose $H \in H^{(r)}(n,p)$ and
\[
E(H) \cap \F \subseteq  \bigsqcup_{i=1}^{t} \prod_{j=1}^r A_j^i,
\]
where `$\sqcup$' denotes the disjoint union. For each $1 \leq i \leq t$,   we assume $A_1^i,A_2^i,\ldots,A_r^i$ form a nontrivial complete $r$-partite $r$-uniform hypergraph.
We fix a constant $K=\tfrac{4}{c} $ and a function  $q(n)=\llll n$.  For each $0 \leq  i \leq  q(n)-1$, we define $f_i=K^i2^{q(n)}$
and
 \[
 \P_i=\left\{ P_i \in \P : f_i \leq  \prod_{j=1}^{r-1}|A_j^i| < f_{i+1} \right\}.
 \]

\begin{lemma} \label{l:lm4}
There is some $0 \leq i  \leq q(n) -1$ such that $|\P_i| \leq \tfrac{t}{q(n)}$.
\end{lemma}

\noindent
{\bf Proof:} Suppose the lemma is not true. As $\P_i$'s are pairwise disjoint, then we have
\[
|\P | \geq \sum_{i=0}^{q(n)-1} |\P_i| > t,
\]
which is a contradiction to the assumption on the size of $\P$. \hfill $\square$

Let  $0 \leq i_0 \leq q(n)-1$ be  the smallest integer satisfying the statement of Lemma \ref{l:lm4}. We consider
\[
 \P'=\left\{ P_i \in \P : \prod_{j=1}^{r-1}|A_j^i|<f_{i_0+1}  \right\}.
\]

\begin{figure}[htbp]
 \centerline{\psfig{figure=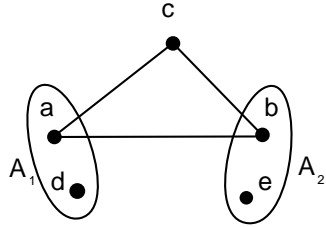, width=0.3\textwidth} }
\caption{An example with $r=3$, $F=\{a,b,c\}$, and $P=\{A_1,A_2\} \in N_{\P',F}(c)$.}
\label{fg:fg2}
\end{figure}

 For an $r$-set $F=\{v_1,v_2,\ldots,v_r\} \in \F$ and each $v_j \in F$, we define
\[
N_{\P',F}(v_j) =\left\{ P_i \in \P' : v_j  \not \in \cup_{s=1}^{r-1} A_s^i, \textrm{ and } |F \cap A_s^i|=1 \textrm{ for each } 1 \leq s \leq r-1 \right\}.
\]
Figure \ref{fg:fg2} is an example for $P \in N_{\P',F}(v_j)$.

Roughly speaking, each $P_i \in N_{\P',F}(v_j)$ could be the prefix of a possible nontrivial complete $r$-partite $r$-uniform hypergraph which contains $F$.

We note that $N_{\P',F}(v_j)$ and $N_{\P',F}(v_k)$ are disjoint if $j \not = k$.
Let $N_{\P'}(F)=\cup_{j=1}^r N_{\P',F}(v_j)$  and $d_{\P'}(F)=|N_{\P'}(F)|= \sum_{j=1}^{r}|N_{\P',F}(v_j)|$.

\begin{lemma} \label{l:lm5}
Assume $|\P' \setminus \P_{i_0}|=xn^{r-1}$ with $x \geq 0.01c$.  Let $\F'=\{F \in \F : d_{\P'}(F) \leq \tfrac{3}{c} x f_{i_0}\}$. We have
\[
|\F'|  \geq \frac{cn^r}{3}.
\]
\end{lemma}

\noindent
{\bf Proof:} We observe that  each $P_i=\{A_1^i,\ldots,A_{r-1}^i\} \in \P'$  contributes one to $d_{\P'}(F)$ for at most $n \prod_{j=1}^{r-1}|A_j^i|$ $r$-sets  $F \in \F$.  Recall the definition of $\P_i$ and Lemma \ref{l:lm4}. For $n$ large enough, we have
\begin{align*}
\sum_{F \in \F} d_{\P'}(F) & \leq \sum_{P_i \in \P'}  n \prod_{j=1}^{r-1} |A_j^i|\\
                          &= \sum_{P_i \in \P' \setminus \P_{i_0}} n \prod_{j=1}^{r-1} |A_j^i| + \sum_{P_i \in \P_{i_0}} n \prod_{j=1}^{r-1} |A_j^i|\\
                          & \leq   nf_{i_0}|\P' \setminus \P_{i_0}|+  nf_{i_0+1}|\P_{i_0}| \\
                          &\leq  xf_{i_0}n^r + \frac{tn f_{i_0+1} }{q(n)} \\
                          & \leq 2x f_{i_0} n^r
\end{align*}
we applied    facts $t \leq n^{r-1}$ and $x \geq 0.01c$ as well as  the definition of $i_0$. We get the following inequality
\[
\frac{3x}{c} f_{i_0} |\F \setminus \F'|  \leq \sum_{F \in \F \setminus \F'} d_{\P'}(F) \leq \sum_{F \in \F} d_{\P'}(F) \leq  2xf_{i_0} n^r .
\]
Clearly, the inequality above implies $|\F \setminus \F'| \leq \tfrac{2cn^r}{3}$. Equivalently, $|\F'| \geq \tfrac{cn^r}{3}$. \hfill $\square$

Before proving a lower bound on the number of uncovered edges, we need one more lemma.

\begin{lemma} \label{l:lm51}
Let $\F'$ be the subfamily of $\F$ given by  Lemma \ref{l:lm5}.  There is a subset  $\W \subseteq \F'$ and a collection of $r$-sets  $\C(F)  \subset \binom{[n]}{r} \setminus \W$   associated with  each $F \in \W$ which satisfy following

\noindent
1:  $|\W|  \geq \tfrac{c^2n^r}{10xf_{i_0}}.$

\noindent
2:  $|\C(F)| \leq \tfrac{3}{c}xf_{i_0}$ for each $F \in \W$.

\noindent
3: For each $F=\{v_1,\ldots,v_r\} \in \W$ and each $1 \leq i \leq r$, if $P=\{A_1,\ldots,A_{r-1}\} \in \P' $ and $P \in N_{\P',F}(v_i)$, then there is  $w \in  A_{r-1} \setminus F$ such that $(F\setminus v_i) \cup w \in \C(F)$.
\end{lemma}

\noindent
{\bf Proof:}   To define $\W$,
 we first give a linear ordering of $r$-sets in $\F'$ and consider the following algorithm.
We will define sets $\F_i$  recursively and build the set $\W$ step by step.
Initially, let  $\F_0=\F'$ and $\W=\emptyset$.

For each $i \geq 1$, if $\F_{i-1} \not = \emptyset$, then let $F_i=\{v_1,v_2,\ldots,v_r\}$ be the first $r$-set in $\F_{i-1}$.
 For each $1 \leq j \leq r$  and each $P=\{A_1,\ldots, A_{r-1}\} \in N_{\P',F_i}(v_j)$,  we note $|F \cap A_s|=1$ for each $1 \leq s \leq r-1$ and $v_j \not \in \cup_{s=1}^{r-1} A_s$ by the definition of $P \in N_{\P', F_i}(v_j)$. Suppose $F_i \cap A_{r-1}=u$.
  We notice $|A_{r-1}| \geq 2$ as $P$ is the prefix of a nontrivial complete $r$-partite $r$-uniform hypergraph.
  If $\left( F_i \setminus u \right) \cup v   \not \in \F_{i-1} \cup \W $ for some $v \in A_{r-1}$,  then  for each $w \in A_{r-1} \setminus v$ , we move   $(F_i \setminus u) \cup w$ from $ \F_{i-1}$  to $\W$ provided $(F_i \setminus u) \cup w \in \F_{i-1}$.
Otherwise,   $\left( F_i \setminus u \right) \cup v   \in \F_{i-1} \cup \W $  for each $v \in A_{r-1}$.  We claim actually $\left( F_i \setminus u \right) \cup v \in \F_{i-1} $ for each $v \in A_{r-1}$. We proceed with the algorithm by assuming this claim.
We choose an arbitrary $w \in A_{r-1} \setminus u$ and delete $\left( F_i \setminus u \right) \cup w $ from $\F_{i-1}$. Moreover, for each $v \in A_{r-1} \setminus w$, we move $\left( F_i \setminus u \right) \cup v $ from $\F_{i-1}$ to $\W$.  We define  $\F_i$ as the resulted subset of $\F_{i-1}$ for each case.  

\begin{figure}[htbp]
 \centerline{\psfig{figure=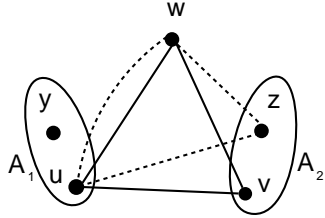, width=0.3\textwidth} }
\caption{An example with $r=3$, $F=\{u,v,w\}$, $P=\{A_1,A_2\} \in N_{\P'}(F)$, and $\{u,w,z\} \in \C(F)$.}
\label{fg:fg3}
\end{figure}

Now we prove the claim. Suppose there is some $z \in A_{r-1} \setminus u$ such that $(F_i \setminus u) \cup z \in W$.  We pick  such a vertex $z$ so that the $r$-set $F'=(F_i \setminus u) \cup z $ is the smallest one in $W$ under the linear ordering.  Suppose $F'$ was added to $W$ at step $j$ with $j<i$.   We examine the moment  that $F'$ was moved to $\W$. If there is some $s \in A_{r-1} \setminus z$ such that $(F' \setminus z) \cup s \not \in \F_{j-1} \cup \W$, then $s \not = u$. Otherwise, $F_i \not \in \F_{i-1}$. We notice $(\F_{i-1} \cup \W) \subseteq (\F_{j-1} \cup \W)$ as $j<i$. Therefore, $(F_i \setminus u) \cup s \not \in \F_{i-1} \cup \W$ and we are in the first case, which is a contradiction.  We obtain $F'$ must satisfy the statement of the claim (because $F'$ is the first one of the form $(F_i \setminus u) \cup z$).
 Thus  $F_i$  was moved to   $\W$ when  we were moving $F'$   to $\W$,  which leads a contradiction.   Figure \ref{fg:fg3} is an illustrative example.

If $\F_{i-1} = \emptyset$, then we halt the process and output $\W$.

Recall the definition of $\F'$, i.e., $d_{\P'}(F) \leq \tfrac{3}{c} x f_{i_0}$ for each $F \in \F'$.  We get that each $F \in \F'$ can make at most  $\tfrac{3}{c}xf_{i_0}$ other $r$-sets in $\F_{i-1}$  deleted from $\F_{i-1}$ if $F$ is added to $\W$ at time $i$. Recall $|\F'| \geq \tfrac{cn^r}{3}$. Thus,
\[
|\W| \geq \frac{|\F'|}{\tfrac{3}{c}xf_{i_0}+1} \geq \frac{c^2n^r}{10xf_{i_0}}.
\]
For each $F \in \W$, we next associate with $F$ a set of $r$-sets $\C(F) \subset \binom{[n]}{r} \setminus \W$. Assume $F=\{v_1,\ldots,v_r\}$.
For each $1 \leq i \leq r$ and each  $\{A_1,\ldots,A_{r-1}\} \in N_{\P',F}(v_i)$,  as the construction of $\W$, there is some $v \in A_{r-1} \setminus  F$ such that $\left( F \setminus v_i \right) \cup v  \not \in \W$. The desired vertex $v$ exists by considering when $F$ is moved to $\W$. If $\left( F \setminus v_i \right) \cup v $ is not an edge, then it excludes the possibility that $F$ get covered by the complete $r$-partite $r$-uniform hypergraph with the prefix $\{A_1,\ldots,A_{r-1}\}$.
 We put the $r$-set  $\left( F \setminus v_i \right) \cup v $ in $\C(F)$.
 For an example, see Figure \ref{fg:fg3}. We will call each $R \in \C(F)$ a {\it certificate} for $F$.
We note that if $R \in \C(F)$, then $|F \cap R|=r-1$ and the symmetric difference $F \bigtriangleup R$  is in  $A_{r-1}$.
We have $|\C(F)| \leq \frac{3}{c}xf_{i_0}$ as the assumption for $|N_{\P'}(F)|$ for each $F \in \F'$.  \hfill $\square$

The next lemma will tell us that with high probability the number of $r$-sets $F \in \W$ such that $F$ is an edge in $H \in H^{(r)}(n,p)$ and $F$ is not contained in any nontrivial complete $r$-partite $r$-uniform hypergraph with the prefix from $\P'$ is large.


\begin{lemma} \label{l:lm6}
 Assume  $1/\llll n \leq p \leq 1/2$, $|\P' \setminus \P_{i_0}|=xn^{r-1}$ with $x \geq 0.01c$,  $H \in H^{(r)}(n,p)$, and Lemma \ref{l:lm1} holds. With probability at most $2\exp(-n^{r-0.92})$,  the number of edges in $E(H) \cap \F $ which is not contained in any complete $r$-partite $r$-uniform hypergraph with the prefix from $\P'$ is less than
\[
\frac{c^2 n^r p(1-p)^{\tfrac{3}{c}xf_{i_0} }}{12 xf_{i_0}}.
\]
\end{lemma}

\noindent
{\bf Proof:}  We will work on the collection of $r$-sets $\W$ given by Lemma \ref{l:lm51}. Let $Y$ be  the number of $r$-sets from $\W$ which is an edge in $H \in H^{(r)}(n,p)$ and is  not contained by any complete $r$-partite $r$-uniform hypergraph with the prefix from $\P'$.  For each $F \in \W$ and  $R \in \C(F)$, we recall that $R$ is a certificate for $F$. We remark that an $r$-set $R$ could be a certificate for more than one $r$-set $F \in \W$.

Let $\C=\cup_{F \in \W} \C(F)$.
For an $r$-set $R \in \C$, if $R \in \C(F)$ for more than $n^{0.45}$ sets $F \in \W$, then we call $R$  a {\it bad certificate}. Let $\C_1$ be the collection of bad certificates. For each $F \in \W$, we set $\C'(F)=\C(F) \setminus \C_1$.  We fix the selection of $\W$,  $\C'(F)$ for each $F \in \W$, and the collection of bad certificates $\C_1$.  We sample all possible edges
 and let $X_F$ be the indicator random variable such that $F$ is an edge in $H \in H^{(r)}(n,p)$ and $R$ is not a edge in $H \in H^{(r)}(n,p)$ for each $R \in \C'(F)$.
  We observe that if $X_F=1$, then $F$ is not covered by any nontrivial complete $r$-partite $r$-uniform hypergraphs with the  prefix   from $\P'$  and   containing  no bad certificate.
 To see this, suppose $F$ is covered by  some $H$ with vertex parts $A_1,\ldots,A_{r-1},A_r$ and $H$ dose not contain any bad certificate. Since the definition of $\C'(F)$,  there is some  $F' \in \C'(F) \cap \prod_{i=1}^{r} A_i $. As the definition of $X_F=1$, we get $F'$ is not an edge. Thus,  $A_1,\ldots, A_r$ do not form a complete $r$-partite $r$-uniform hypergraph, which  is a contradiction.

 We define $X=\sum_{F \in \W} X_F$.  Applying Lemma \ref{l:lm3} with $\F=\W$,  $k(n)=\tfrac{3}{c}xf_{i_0}$, and $l(n)=n^{0.45}$, we obtain
with probability at least $1-2\exp(-n^{r-0.92})$, we have
\[
X \geq \frac{c^2n^r p(1-p)^{\tfrac{3}{c}xf_{i_0}}}{11xf_{i_0}}.
\]
We note $n^{r-0.01}$ is a lower term as the definition of $f_{i_0}$ and the assumption for $p$.
We use $\F''$ to denote those  $r$-sets $F \in \W$ such that  $X_F=1$.
The argument above gives that with probability at least $1-2\exp(-n^{r-0.92})$, we have
\[
|\F''|  \geq \frac{c^2n^r p(1-p)^{\tfrac{3}{c}xf_{i_0}}}{11xf_{i_0}}.
\]
Let us condition on this.

 We note that an edge in $\F''$ could be covered by some complete $r$-partite $r$-uniform hypergraph which contains a bad certificate. We next prove an upper bound on the number of such  edges. This upper bound  works for all samplings of edges.

Let $A_1,\ldots,A_r$ be the vertex parts of such a complete $r$-partite $r$-uniform hypergraph $H$. Suppose $\{A_1,\ldots,A_{r-1}\} \in \P'$.
We define
\[
A_r'=\{v_r \in A_r:  \textrm{ there are } v_1 \in A_1, \ldots,v_{r-1} \in A_{r-1} \textrm{ such that } \{v_1,\ldots,v_r\} \in \F''  \}.
\]
The number of edges from $\F''$ covered by $H$ is at most $|A_i'| \prod_{i=1}^{r-1}|A_i|$. We next relate the number of bad certificates contained in $H$ to the size of $A_r'$.

For each $w \in A_r'$, by the definition of $A_r'$, there is  some $F=\{v_1,\ldots,v_{r-1},w\} \in A_1 \times \cdots \times A_{r-1} \times \{w\}$ such that  $F \in \F''$.
We observe $\{A_1,\ldots,A_{r-1}\} \in N_{\P'}(F)$.
Let $\{v_1,\ldots,v_{r-2},z,w\}$ be the certificate of $F$ associated with $\{A_1,\ldots,A_{r-1}\}$, where
 $z \in A_{r-1}$. We notice $\{v_1,\ldots,v_{r-2},z,w\}$ must be a bad certificate.
  Otherwise, as $F \in \F''$, we get $\{v_1,\ldots,v_{r-2},z,w\}$ is a nonedge. Then $A_1,\ldots,A_r$ do not form a complete $r$-partite $r$-uniform hypergraph which is a contradiction.
      Therefore, each $w \in A_r'$ gives  at least one bad certificate from $A_1 \times \cdots \times A_{r-1} \times \{w\}$ and these bad certificates are distinct for different $w$ from $A_r'$. We obtain that the number of bad certificate in $H$ is at least $|A_r'|$.

We divide those hypergraphs which contains a bad certificate into two subsets $\h_1$ and $\h_2$, where $\h_1=\{H: |A_r'| \leq n^{0.9}\}$ and $\h_2=\{H: |A_r'| \geq n^{0.9}\}$. We note that each $H \in \h_2$ contains at least $n^{0.9}$ bad certificates as the analysis above.
We next prove absolute upper bounds for the number of edges from $\F''$ which are covered by $\h_1$ and $\h_2$.
We observe that each $H \in \h_1$ can cover at most $ |A_r'| \prod_{j=1}^{r-1} |A_j| \leq (r+1) n^{0.9}\log n $ edges from $\F''$ as we assume Lemma \ref{l:lm1} holds. There are at most $t<n^{r-1}$ of them as assumptions in Theorem \ref{t:submain}. Therefore, $\h_1$ covers at most $(r+1)n^{r-0.1}\log n$ edges from $\F''$.

We need an upper bound for the number of bad certificates in total. We consider pairs $(F,R)$ such that $F \in \W$ and $R \in \C(F)$. As
$|\C(F)| \leq \tfrac{3}{c} x f_{i_0}$ for each $F \in \W$, the number of such pairs is less than
\[
|\W| \frac{3}{c} x f_{i_0}< \frac{3}{c}xf_{i_0} n^r < n^{r} \log n,
\]
here we used the fact $|\W| < n^r$ and the definition of  $f_{i_0}$.
 As the definition of a bad certificate, a simple double counting method yields that the number of bad certificates is at most $n^{r-0.55} \log n$. Since each bad certificate (viewed as an edge) is contained in at most one $H \in \h_2$ (we are considering the partition of edges) and  each $H \in \h_2$ contains at least $n^{0.9}$ bad certificates, we have $|\h_2| \leq n^{r-1.45} \log n$. The number of edges contained in each $H \in \h_2$ has an absolute upper bound $(r+1)n\log n$. Therefore, the number of edges from $\F''$ which are covered by $\h_2$ is at most $(r+1)n^{r-0.45} \log^2 n$.

Thus, those complete $r$-partite $r$-uniform hypergraphs containing a  bad certificate cover at most $(r+1)n^{r-0.1}\log n+(r+1)n^{r-0.45} \log^2 n$ edges from $\F''$.  Therefore,  we have
\begin{align*}
Y & \geq \frac{c^2n^r p(1-p)^{\tfrac{3}{c}xf_{i_0}}}{11xf_{i_0}}-(r+1)n^{r-0.1}\log n-(r+1)n^{r-0.45} \log^2 n\\
   & > \frac{c^2n^r p(1-p)^{\tfrac{3}{c}xf_{i_0}}}  {12xf_{i_0}}
\end{align*}
the proof of this lemma is complete.
\hfill $\square$

We are now ready to prove Theorem \ref{t:submain}.

\noindent
{\bf Proof of Theorem \ref{t:submain}:} To simplify the notation, we define the following prefix sets:
\begin{align*}
\Q_1=\P' \setminus \P_{i_0} &=\left\{P_i \in \P: \prod_{j=1}^{r-1} |A_j^i| < f_{i_0} \right\} \\
\Q_2=(\P\setminus \P') \cup \P_{i_0} &=\left\{P_i \in \P: \prod_{j=1}^{r-1}  |A_j^i| \geq f_{i_0}\right\} \\
\Q_3=\P \setminus \P' &=\left\{P_i \in \P: \prod_{j=1}^{r-1} |A_j^i|  \geq f_{i_0+1} \right\} \\
\end{align*}
Let $c_1(n)=2, c_2(n)=f_{i_0},$ and $c_3(n)=f_{i_0+1}$. For $H \in H^{(r)}(n,p)$ and each $i \in \{1,2,3\}$, let ${\cal Z}_i$ be the event that $g(H,\Q_i) \leq |\Q_i|c_i(n)p^{c_i(n)}n+2n^{r-0.3}$. Lemma \ref{l:lm2} implies that with probability at least $1-6\exp(-n^{r-0.8})$ all events ${\cal Z}_1,{\cal Z}_2,{\cal Z}_3$ hold simultaneously. We condition on these three events. We note that for each $i \in \{1,2,3\}$,  the number of edges from $\F$ which is covered by complete $r$-partite $r$-uniform hypergraphs with the prefix from $\Q_i$ is bounded above by the function $g(H,\Q_i)$.

 We proceed to prove  $|\Q_1| \geq 0.01c$.
Suppose not.  Because  the event ${\cal Z}_1$  occurs,  the number of edges from $\F$ covered by those complete $r$-partite $r$-uniform hypergraphs  with the prefix $P \in \Q_1$ is at most $(2+o(1))p^2 n|\Q_1| \leq \left(0.02cp^2+o(1) \right) n^r $, here $2n^{r-0.2}$ is a lower term as we assume $p \geq 1/\llll n$.
 A simple application of Theorem \ref{t:chernoff} yields that with probability at least $1-\exp(-cpn^r/8)$  the number of $r$-sets in $\F$   being an edge in   $H \in H^{(r)}(n,p)$ is at least  $\tfrac{cpn^r}{2}$. Therefore, the number of edges covered by those  complete $r$-partite $r$-uniform hypergraphs with the prefix $\Q_2$  is at least $\tfrac{cp n^r}{4} $. As the event ${\cal Z}_2$, we get
 \[
|\Q_2| \geq \frac{(\tfrac{pc}{4}+o(1))n^r}{ f_{i_0} p^{f_{i_0}} n} > n^{r-1}
\]
 when $n$ is large enough. This is a contradiction to the assumption $|\P| \leq n^{r-1}$. Therefore, as long as events ${\cal Z}_1$ and ${\cal Z}_2$  as well as the lower bound for the number of edges from $\F$ hold,  we have $|\Q_1| \geq 0.01c$ which is one of the assumptions in Lemma \ref{l:lm6}.

Recall Lemma \ref{l:lm6}. Those uncovered edges given by Lemma \ref{l:lm6} must be covered by complete $r$-partite $r$-uniform hypergraphs  with the prefix from $\Q_3$. As  the event ${\cal Z}_3$, we get
\[
|\Q_3| \geq  \frac{c^2 n^{r-1} p(1-p)^{\tfrac{3}{c}xf_{i_0}}}{13 x f_{i_0} f_{i_0+1} p^{f_{i_0+1}}},
\]
we note $n^{r-0.3}$ is a lower order term. Recall $1/\llll n \leq p \leq 1/2$ and $f_i=K^i 2^{q(n)}$.
 We get
\begin{align}
|\P|=|\P'|+|\Q_3| & \geq xn^{r-1}+\frac{c^2 n^{r-1} p(1-p)^{\tfrac{3}{c}xf_{i_0}}}{13 x f_{i_0} f_{i_0+1} p^{f_{i_0+1}}} \nonumber\\
&\geq  \frac{c^2 n^{r-1} p(1-p)^{\tfrac{3}{c}xf_{i_0}}}{13 x f_{i_0} f_{i_0+1} p^{f_{i_0+1}}} \nonumber \\
& \geq \frac{c^2 n^{r-1} p 2^{f_{i_0+1}-\tfrac{3}{c}xf_{i_0}}}{13xf_{i_0}f_{i_0+1}} \label{lb1}\\
& \geq \frac{c^2 n^{r-1} p 2^{f_{i_0+1}-\tfrac{3}{c}f_{i_0}}}{13f_{i_0}f_{i_0+1}} \label{lb2}\\
&= \frac{c^2 n^{r-1} p 2^{K^{i_0+1}2^{q(n)}-\tfrac{3}{c}K^{i_0}2^{q(n)}}}{13K^{2i_0+1}2^{2q(n)}} \nonumber\\
&=\frac{c^2 n^{r-1} p 2^{\tfrac{1}{c}K^{i_0}2^{q(n)}}}{13K^{2i_0+1}2^{2q(n)}} \label{lb3}\\
&> n^{r-1} \nonumber
\end{align}
when $n$ is large enough. We used $p \leq 1/2$ to get inequality \eqref{lb1}, $x \leq 1$ to get inequality \eqref{lb2}, and $K=\tfrac{4}{c}$ to get inequality \eqref{lb3}.

Therefore, as long as events ${\cal Z}_1, {\cal Z}_2, {\cal Z}_3$ occur,   the lower bound for the number of edges in $\F \cap E(H)$ holds, and Lemma \ref{l:lm6} holds, we get a contradiction. With probability at most $2\exp(-n^{r-0.92})+6\exp(-n^{r-0.8})+\exp(-cpn^r/8) \leq 3\exp(-n^{r-0.92})$,
one of them does not hold, this completes the proof of the theorem.
\hfill $\square$

\section{Proof of Theorem \ref{main}}
Before we prove the main theorem, we need to show  an upper bound on the number of choices for the prefix set $\P$.
\begin{lemma} \label{l:lm10}
Suppose $\P=\{P_1,\ldots,P_q\}$, where  $P_i=\{A_1^i,\ldots,A_{r-1}^i\}$  and $1 \leq \prod_{j=1}^{r-1} |A_j^i| < (r+1)\log n$ for each $1 \leq i \leq q$.
The number of choices for $\P$ with $|\P| \leq  n^{r-1}$ is bounded from above by $n^{(r+3)n^{r-1}\log n}$ when $n$ is large enough.
\end{lemma}

\noindent
{\bf Proof:} We shall show the desired upper bound step by step. We have at most $n^{r-1}$ choices for the size of $\P$. First, we fix the size of $\P$. We will  establish an absolute  upper bound on the number of  choices for each element $P_i$ of $\P$.  For each $P_{i}=\{A_1^i,\ldots,A_{r-1}^i\} \in \P$, we have $t_i=|\cup_{j=1}^{r-1} A_j^i|  \leq (r+1)\log n+s <  (r+2)\log n$ as $\prod_{j=1}^{r-1} |A_j^i| \leq (r+1)\log n$. Therefore,
$\cup_{j=1}^{r-1} A_j^i \in   \binom{[n]}{\leq (r+2) \log n}$, which implies that the number of choices for  $\cup_{j=1}^{r-1} A_j^i$ is at most $n^{(r+2)\log n}$. We fix  the selection of $\cup_{j=1}^{r-1} A_j^i$ and wish to  partition it into $r-1$ disjoint parts $A_j^i$. Let $a_j=|A_j^i|$ for $1 \leq j \leq r-1$. Then we have $a_1+\ldots+a_{r-1}=t_i$. The number of choice for the size of $a_1,\ldots,a_{j-1}$ equals the number of solutions to the equation $a_1+\ldots+a_{j-1}=t_i$.  Since $a_j \geq 1$,  we have at most $\binom{t_i}{r-1}$ choices for $a_1,\ldots,a_{j-1}$
, which can be bounded from above by $((r+2)\log n)^{r-1}$ as $t_i \leq (r+2)\log n$. If we fix the size of each $A_j^i$, then the number of ways to partition  $\cup_{j=1}^{r-1} A_j^i$ into $A_1^i,\ldots,A_{r-1}^i$ equals  $\binom{t_i}{a_1,\ldots,a_{j-1}}$, which is at most $t_i! \leq ((r+2)\log n)^{(r+2)\log n}$. Therefore, the number of choices for $P_i$ is at most
\[
n^{(r+2)\log n}\left((r+2)\log n  \right)^{(r+2)\log n+r-1}
\]
Recall the assumption  $|\P| \leq n^{r-1}$. We get that the number of choices for $\P$ is at most
\[
n^{r-1} \left(n^{(r+2)\log n}\left((r+2)\log n  \right)^{(r+2)\log n+r-1}  \right)^{|\P|}     <n^{(r+3)n^{r-1}\log n},
\]
provide $n$ is sufficiently large.
\hfill $\square$

We are ready to prove  Theorem \ref{main}.

\noindent
{\bf Proof of the upper bound:}  We shall exhibit an explicit decomposition of each $r$-uniform hypergraph with $n$ vertices using at most $(1-\pi(K_r^{(r-1)})+o(1))\binom{n}{r-1}$ trivial complete $r$-partite $r$-uniform hypergraphs.

 For each $r \geq 3$,  let $G=([n],E)$ be an $(r-1)$-uniform hypergraph which has $\textrm{ex}(n,K_{r}^{(r-1)})$ edges and does not contain  $K^{(r-1)}_{r}$ as a subhypergraph. Obviously, $G$ is well-defined.
 Let $G'$ be the complement of $G$. Therefore, $E(G')=\binom{[n]}{r-1} \setminus E(G)$. We observe that an independent set of size $r$ in $G'$ will be a $K_{r}^{(r-1)}$ in $G$. As $G$  does not contain $K_{r}^{(r-1)}$, we get each $F \in \binom{[n]}{r}$ contains at least one edge of $G'$.

  Suppose $q=|E(G')|$ and we list edges in $G'$ as $e_1,\ldots,e_q$. For each $r$-uniform hypergraph $H$ with $n$ vertices, we will show that $H$ can be decomposed into at most $q$ complete $r$-partite $r$-uniform hypergraphs as follows.

  Let $H_0=H$ and we will define a sequence of  complete $r$-partite $r$-uniform hypergraphs recursively.  For each $1 \leq i \leq q$, we assume the edge $e_i$ in $G'$ is
 $\{v_1,v_2,\ldots,v_{r-1}\}$. The key observation is the following. For an edge $F \in E(H)$, if $F$ is contained in a trivial complete $r$-partite $r$-uniform hypergraphs with vertex parts $\{v_1\} \times \cdots \times \{v_{r_1}\} \times V_r$, then the set $\{v_1,\ldots,v_{r-1}\}$ must be a subset of  $F$. We define $\F_i=\{F \in E(H_{i-1}): e_i \subset F\}$ and $A_{r}=\cup_{F \in \F_i} F \setminus e_i$. If $A_{r}\not = \emptyset$, then the $i$-th complete $r$-partite $r$-uniform hypergraphs $H_i'$ will have vertex parts $\{v_1\},\ldots, \{v_{r_1}\}, A_r$. If the set $\F_i$ is empty, then we do not define $H_i'$. We set $E(H_i)=E(H_{i-1}) \setminus E(H_i')$ for each $1 \leq i \leq q-1$.

We note that each $F \in E(H)$ contains at least one   of $e_i$. The definition of $H_i$'s ensures that each edge in $H$ is in  exactly one of these trivial complete $r$-partite $r$-uniform hypergraphs. Clearly, for sufficiently large $n$, we have $q=(1-\pi(K^{(r-1)}_{r})+o(1)) \binom{n}{r-1}$. Since the decomposition above applies  to all $H$, it also works for the random hypergraph $H \in H^{(r)}(n,p)$.

 \vspace{0.1cm}

 \noindent
 {\bf Proof of the lower bound:} We assume Lemma \ref{l:lm1} holds.  Thus  each complete $r$-partite $r$-uniform hypergraph with vertex parts $A_1,A_2,\ldots,A_{r}$ satisfying $\prod_{i=1}^{r-1} |A_r| < (r+1) \log n$ provided $|A_1| \leq |A_2| \leq \ldots \leq |A_{r}|$.
 For any fixed small positive constant $\epsilon$, we shall show that the probability $f(H) \leq (1-\pi(K_r^{(r-1)})-\epsilon) \binom{n}{r-1}$ is small, where $H \in H^{(r)}(n,p)$.

  For a fixed prefix set $\P=\{P_1,\ldots,P_t\}$, where
$P_i=\{A_1^i,\ldots,A_{r-1}^i\}$ and $\prod_{j=1}^{r-1} |A_j^i| < (r+1) \log n$ for each $1 \leq i \leq t$ and $t \leq (1-\pi(K_{r}^{(r-1)})-\epsilon)\binom{n}{r-1}$. Let $\X$ denote the event that there are $t$ sets $A_r^1,\ldots,A_r^t$ such that
\[
E(H)=\bigsqcup_{i=1}^t \prod_{j=1}^{r} A_j^i,
\]
provided $H \in H^{(r)}(n,p)$. Here $A_1^i,\ldots,A_{r}^i$ form a complete $r$-partite $r$-uniform hypergraph for each $1 \leq i \leq t$.
 We assume the first $s$ of them are trivial complete $r$-partite $r$-uniform hypergraphs, i.e., $|A_1^i|=\cdots=|A_{r-1}^i|=1$ for each $1 \leq i \leq s$.

  As we did for proving the upper bound, we define an $(r-1)$-uniform hypergraph $G$ such that $V(G)=[n]$ and $E(G)= \binom{[n]}{r-1} \setminus (\cup_{i=1}^s \prod_{j=1}^{r-1} A_j^i)$.  We note $|A_j^i|=1$ for each $1 \leq i \leq s$ and $1 \leq j \leq r-1$. We get $|E(G)| \geq (\pi(K_{r}^{(r-1)})+\epsilon) \binom{n}{r-1}$. By the supersaturation result for hypergraphs (see Theorem 1 in \cite{es}), we get that  there are at least $c(\epsilon) n^{r}$ copies of $K_{r}^{(r-1)}$ in $G$. Let $G'$ be the complement of $G$ and $\F$ be the collection of independent sets with size $r$ in $G'$. We have $|\F| \geq c(\epsilon)n^{r}$.  We observe that if $H \in H^{(r)}(n,p)$, then edges in $\F \cap E(H)$ must be covered by  those nontrivial complete $r$-partite $r$-uniform hypergraphs in the partition.
   Let $\cal Y$ be the event that each $F \in \F \cap E(H)$ is contained in exactly one of  the last $t-s$ nontrivial complete $r$-partite $r$-uniform hypegraphs.
  We have two cases depending on the range of the probability $p$.
\begin{description}
\item[Case 1:]  $1/\llll n \leq p \leq 1/2$.  Applying Theorem \ref{t:submain} with $\P'=\{P_{s+1},\ldots,P_{t}\}$, we get  that $\cal Y$ holds with the probability at most $2\exp(-n^{r-0.92})$. This implies that the event $\X$ occurs with probability at most $2\exp(-n^{r-0.92})$. By Lemma \ref{l:lm10}, there are at most $n^{(r+3)n^{r-1} \log n}$ choices for $\P$ satisfying the desired properties. Applying the union bound, we get that the probability $f(H) \leq (1-\pi(K_r^{(r-1)})-\epsilon) \binom{n}{r-1}$ is at most
     $2\exp(-n^{r-0.92}) n^{(r+3)n^{r-1} \log n}< \exp(-n^{r-0.94})$ for any positive $\epsilon$.

\item[Case 2:] $(\log n)^{2.001}/n \leq p \leq 1/\llll n$.   We observe that the set $\F$ is determined by the prefix set $\P$. Therefore, Lemma \ref{l:lm10} also gives an upper bound on the number of possible choices of $\F$.   A simple application of Theorem \ref{t:chernoff}  yields that   with high probability $|\F \cap E(H)| \geq \tfrac{pcn^r}{4}$ for all $\F$ with $|\F| \geq n^r/\log \log n$.  As edges in $\F \cap E(H)$ must be covered by the last $t-s$ nontrivial complete $r$-partite $r$-uniform hypergraphs. Since $t-s \leq \binom{n}{r-1}$,   Lemma \ref{l:lm9} tells us that the event $\cal Y$ occurs with probability at most $2^{-0.05pc\log(1/p)n^r}$. This also implies that the event $\X$ occurs with probability at most $2^{-0.05pc\log(1/p)n^r}$. By Lemma \ref{l:lm9} and the union bound, we get the probability $f(H) \leq (1-\pi(K_r^{(r-1)})-\epsilon) \binom{n}{r-1}$ is at most
 $2^{-0.05pc\log(1/p)n^r} n^{(r+3)n^{r-1} \log n} \leq 2^{-0.04pc\log(1/p)n^r}$ as $np \geq (\log n)^{2.001}$ and $c$ is a constant.

\end{description}
The proof of the theorem is finished.
  \hfill $\square$

\section{Concluding remarks}
In this paper, we studied the problem of partitioning the edge set of a random $r$-uniform hypergraph into edge sets of complete $r$-partite $r$-uniform hypergraphs.  We were able to show
if $(\log n)^{2.001}/n \leq p \leq 1/2$ and $H \in H^{(r)}(n,p)$, then with high probability $f(H)=(1-\pi(K^{(r-1)}_r)+o(1))\binom{n}{r-1}$.  For the case of $r=2$, results from \cite{alon2} and \cite{cp} assert that if  $p$ is a constant, $p \leq 1/2$,  and $G \in G(n,p)$, then with high probability $n-o((\log n)^{3+\epsilon}) \leq  f(G) \leq  (2+o(1))\log n$ for any  positive constant $\epsilon$.  For sparse random graphs, Alon \cite{alon2} determined the order of magnitude of  the second term of $f(G)$.  However, we do not have any information on the second term of $f(H)$ for $r \geq 3$.  This leads the following question.

\vspace{.1in}

\noindent
{\bf  Problem 1:}  Determine the order of the second term of $f(H)$ for $H \in H^{(r}(n,p)$ and $r \geq 3$.

We note that we were only able to determine the leading coefficient of $f(H)$ for   $p \geq (\log n)^{2.001}/n$ and $H \in H^{(r}(n,p)$. A natural question is to prove similar results for  other range of the  probability $p$.

We recall that for a graph $G$, the {\it strong bipartition number} $\bp'(G)$  of $G$ is the minimum number of nontrivial complete bipartite subgraphs (which are not stars) of $G$ such that each edge of $G$ is in exactly one of them. This parameter was introduced by Chung and the author in \cite{cp} when they were studying the bipartition number of random graphs. In particular, they proved that   if  $p$ is a constant, $p \leq 1/2$,  and $G \in G(n,p)$, then  $\bp'(G) \geq 1.0001 n$ with high probability.   For sparse random graphs, Alon \cite{alon2} proved a better lower bound. Namely, he showed with high probability $\bp'(G) \geq 2n$ if $G \in G(n,p)$.  We remark here that our methods for proving Theorem \ref{t:submain} implicitly yield the following theorem.
\begin{theorem}
If  $p$ is a constant, $p \leq 1/2$,  and $G \in G(n,p)$, then with high probability
\[
\frac{\bp'(G)}{n} \to \infty  \textrm{  as  } n \to \infty.
\]
\end{theorem}

\end{document}